\documentclass{amsart}
\usepackage{amsmath}
\usepackage{amsfonts}

\newtheorem{theorem}{Theorem}
\theoremstyle{plain}

\newtheorem{definition}{Definition}

\newtheorem{proposition}{Proposition}

\numberwithin{equation}{section}

\begin{document}
\title[Legendre Trajectories]{Legendre Trajectories of Trans-$S$-Manifolds}
\author{\c{S}aban G\"{u}ven\c{c}}
\address[\c{S}. G\"{u}ven\c{c}]{Balikesir University, Department of
Mathematics\\
Balikesir, TURKEY}
\email[\c{S}. G\"{u}ven\c{c}]{sguvenc@balikesir.edu.tr}
\urladdr{http://sguvenc.baun.edu.tr}

\begin{abstract}
In this paper, we consider Legendre trajectories of trans-$S$-manifolds. We
obtain curvature characterizations of these curves and give a classification
theorem. We also investigate Legendre curves whose Frenet frame fields are
linearly dependent with certain combination of characteristic vector fields
of the trans-$S$-manifold. \medskip
\end{abstract}

\maketitle

\section{Introduction}

Let $(M,g)$ be a Riemannian manifold, $F$ a closed $2$-form and let us
denote the Lorentz force on $M$ by $\Phi $, which is a $(1,1)$-type tensor
field.\ If $F$ is associated by the relation%
\begin{equation}
g(\Phi X,Y)=F(X,Y),\text{ \ }\forall X,Y\in \chi (M),  \label{F}
\end{equation}%
then it is called a\textit{\ magnetic field} (\cite{Adachi-1996}, \ \cite%
{BRCF} and \cite{Comtet-1987}). Let $\nabla $ be the Riemannian connection
associated to the Riemannian metric $g$ and $\gamma :I\rightarrow M$ a
smooth curve. If $\gamma $ satisfies the Lorentz equation%
\begin{equation}
\nabla _{\gamma ^{\prime }(t)}\gamma ^{\prime }(t)=\Phi (\gamma ^{\prime
}(t)),  \label{Lorentz eq}
\end{equation}%
then it is called a \textit{magnetic curve}  for the magnetic field $F$. The
Lorentz equation is a generalization of the equation for geodesics. Magnetic
curves have constant speed. If the speed of the magnetic curve $\gamma $ is
equal to $1$, then it is called a \textit{normal magnetic curve }\cite%
{DIMN-2015}. For extensive information about almost contact metric manifolds
and Sasakian manifolds, we refer to Blair's book \cite{Blair}.

Let $\gamma (s)$ be a Frenet curve parametrized by the arc-length parameter $%
t$ in an almost contact metric manifold $M$. The function $\theta (t)$
defined by $cos[\theta (t)]=g(T(t),\xi )$ is called \textit{the contact
angle function}. A curve $\gamma $ is called a \textit{slant curve} if its
contact angle is a constant \cite{CIL}. If a slant curve is with contact
angle $\frac{\pi }{2}$, then it is called a \textit{Legendre curve }\cite%
{Blair}. Likewise, Cihan \"{O}zg\"{u}r and the present author defined
Legendre curves of $S-$manifolds in \cite{OG-2014}. A curve $\gamma :I$ $%
\rightarrow M=(M^{2n+s},f,\xi _{i},\eta _{i},g)$ is called a Legendre curve
if $\eta _{i}(T)=0$, for every $i=1,...s$, where $T$ is the tangent vector
field of $\gamma $. This definition can be used in trans-$S$-manifolds.

Let $\gamma $ be a curve in an almost contact metric manifold $(M,\varphi
,\xi ,\eta ,g)$. In \cite{LSL}, Lee, Suh and Lee introduced the notions of $%
C $-parallel and $C$-proper curves in the tangent and normal bundles. A
curve $\gamma $ in an almost contact metric manifold $(M,\varphi ,\xi ,\eta
,g)$ is defined to be $C$\textit{-parallel }if $\nabla _{T}H=\lambda \xi $, $%
C$\textit{-proper }if $\Delta H=\lambda \xi $, $C$\textit{-parallel in the
normal bundle }if $\nabla _{T}^{\perp }H=\lambda \xi $, $C$\textit{-proper
in the normal bundle} if $\Delta ^{\perp }H=\lambda \xi $, where $T$ is the
unit tangent vector field of $\gamma $, $H$ is the mean curvature vector
field, $\Delta $ is the Laplacian, $\lambda $ is a non-zero differentiable
function along the curve $\gamma ,$ $\nabla ^{\perp }$ and $\Delta ^{\perp }$
denote the normal connection and Laplacian in the normal bundle,
respectively \cite{LSL}. The present author and Cihan \"{O}zg\"{u}r
generalized this definition for $S$-manifolds in \cite{GO-2019}. In the
present study, this definition will be used in trans-$S$-manifolds as well.

An almost contact metric manifold $M$ is called a \textit{trans-Sasakian
manifold }\cite{Oubina}\textit{\ }if there exist two functions $\alpha $ and 
$\beta $ on $M$ such that%
\begin{equation}
(\nabla _{X}\varphi )Y=\alpha \lbrack g(X,Y)\xi -\eta (Y)X]+\beta \lbrack
g(\varphi X,Y)\xi -\eta (Y)\varphi X],  \label{trans}
\end{equation}%
for any vector fields $X,Y$ on $M$. $C$-parallel and $C$-proper slant curves
of trans-Sasakian manifolds were studied in \cite{OG-2014-2}.

\section{Preliminaries}

Firstly, let us recall framed $f$-manifolds. Let ($M,g$) be a ($2n+s$%
)-dimensional Riemann manifold. It is called \textit{framed metric }$\mathit{%
f}$\textit{-manifold} with a \textit{framed metric }$\mathit{f}$\textit{%
-structure} $(f,\xi _{i},\eta _{i},g)$, $\alpha \in \left\{ 1,...,s\right\}
, $ if it satisfies the following equations: 
\begin{equation}
\begin{array}{cccc}
\varphi ^{2}=-I+\overset{s}{\underset{\alpha =1}{\sum }}\eta ^{\alpha
}\otimes \xi _{\alpha }, & \eta _{i}(\xi _{j})=\delta _{ij}, & f\left( \xi
_{i}\right) =0, & \eta _{i}\circ f=0%
\end{array}%
\end{equation}%
\begin{equation}
g(fX,fY)=g(X,Y)-\overset{s}{\underset{i=1}{\sum }}\eta _{i}(X)\eta _{i}(Y),
\end{equation}%
\begin{equation}
\eta ^{\alpha }(X)=g(X,\xi ).  \label{eq3}
\end{equation}%
Here, $f$ is a ($1,1$) tensor field\textit{\ }of rank $2n$; $\xi
_{1},...,\xi _{s}$ are vector fields; $\eta _{1},...,\eta _{s}$ are $1$%
-forms and $g$ is a Riemannian metric on $M$; $X,Y\in \chi (M)$ and $i,j\in
\left\{ 1,...,s\right\} $ \cite{Nak-1966}. $(f,\xi _{i},\eta _{i},g)$ is
called $S$\textit{-structure},\textit{\ }when the Nijenhuis tensor of $%
\varphi $ is equal to $-2d\eta ^{\alpha }\otimes \xi _{\alpha }$, for all $%
\alpha \in \left\{ 1,...,s\right\} $ \cite{Blair}.

Secondly, the concept of trans-$S$-manifolds is as follows:

A $(2n+s)-$dimensional metric $f-$manifold M is called an \textit{almost
trans-}$S$\textit{-manifold} if it satisfies%
\begin{equation}
\left( \nabla _{X}f\right) Y=\sum_{i=1}^{s}\left[ 
\begin{array}{c}
\alpha _{i}\left\{ g\left( fX,fY\right) \xi _{i}+\eta (Y)f^{2}X\right\} \\ 
+\beta _{i}\left\{ g\left( fX,Y\right) \xi _{i}-\eta _{i}(Y)fX\right\}%
\end{array}%
\right] ,  \label{3.1}
\end{equation}%
where $\alpha _{i},\beta _{i}$ $\left( i=1,...,s\right) $ are smooth
functions and $X,Y\in \chi (M)$ \cite{AFP-2016}. If M is normal, then it is
called a \textit{trans-}$S$\textit{-manifold}. If $s=1$, a trans-$S$%
-manifold becomes a \textit{trans-Sasakian manifold}. In trans-Sasakian
case, the above condition implies normality. But, for $s\geq 2$, this
statement is no longer valid \cite{AFP-2016}. Since%
\begin{equation}
\left[ f,f\right] \left( X,Y\right) +2\sum_{i=1}^{s}d\eta _{i}\left(
X,Y\right) \xi _{i}=\sum_{i,j=1}^{s}\left[ \eta _{j}\left( \nabla _{X}\xi
_{i}\right) \eta _{j}\left( Y\right) -\eta _{j}\left( \nabla _{Y}\xi
_{i}\right) \eta _{j}\left( X\right) \right] \xi _{i},  \label{3.2}
\end{equation}%
and $\left\{ \xi _{i}\right\} _{i=1}^{s}$ is $g-$orthonormal, it is found
that%
\begin{equation*}
\sum_{j=1}^{s}\left[ \eta _{j}\left( \nabla _{X}\xi _{i}\right) \eta
_{j}\left( Y\right) -\eta _{j}\left( \nabla _{Y}\xi _{i}\right) \eta
_{j}\left( X\right) \right] =0
\end{equation*}%
for all $i=1,...,s$. After calculations, one obtains%
\begin{equation}
\nabla _{X}\xi _{i}=-\alpha _{i}fX-\beta _{i}f^{2}X,  \label{3.4}
\end{equation}%
for $i=1,...,s$ \cite{AFP-2016}.

The notion of a Frenet curve is well-known as below:

Let us consider a unit-speed curve $\gamma :I\rightarrow M$ in an $n$%
-dimensional Riemannian manifold $(M,g)$. If there exists orthonormal vector
fields $E_{1},E_{2},...,E_{r}$ along $\gamma $ satisfying 
\begin{eqnarray}
E_{1} &=&\gamma ^{\prime }=T,  \notag \\
\nabla _{T}E_{1} &=&\kappa _{1}E_{2},  \notag \\
\nabla _{T}E_{2} &=&-\kappa _{1}E_{1}+\kappa _{2}E_{3},
\label{Frenetequations} \\
&&...  \notag \\
\nabla _{T}E_{r} &=&-\kappa _{r-1}E_{r-1},  \notag
\end{eqnarray}%
then $\gamma $ is called a\textit{\ Frenet curve of osculating order }$r$,
where $\kappa _{1},...,\kappa _{r-1}$ are positive functions on $I$ and$\
1\leq r\leq n.$

A Frenet curve of osculating order $1$ is a called \textit{geodesic}. A
Frenet curve of osculating order $2$ is a \textit{circle} if $\kappa _{1}$
is a non-zero positive constant. A Frenet curve of osculating order $r\geq 3$
is called a \textit{helix of order }$r$, when $\kappa _{1},...,\kappa _{r-1}$
are non-zero positive constants; a helix of order $3$ is simply called a 
\textit{helix}.

Finally, we can define Legendre curves in trans-$S$-manifolds like:

\begin{definition}
Let $M=(M^{2n+s},f,\xi _{i},\eta _{i},g)$ be a trans-$S$-manifold. Consider
a unit-speed smooth curve $\gamma :I\rightarrow M$ and its unit tangential
vector field $T=\gamma ^{\prime }$. If $\eta _{i}(T)=0$ for all $i=1,2,...,s$%
, then it is called a \textit{Legendre curve}.
\end{definition}

Here are the direct results from the definition:%
\begin{equation*}
f^{2}T=-T,
\end{equation*}%
\begin{equation}
\kappa _{1}\eta _{i}\left( E_{2}\right) +\beta _{i}=0,  \label{kappa1etae2}
\end{equation}%
\begin{equation*}
\left( \nabla _{T}f\right) T=\sum_{i=1}^{s}\alpha _{i}\xi _{i},
\end{equation*}%
which gives us%
\begin{equation}
\nabla _{T}fT=\sum_{i=1}^{s}\alpha _{i}\xi _{i}+\kappa _{1}fE_{2}
\label{nablafT}
\end{equation}

Let us recall what a magnetic curve is and what we mean by trajectory:

Let $M^{2n+s}=(M^{2n+s},f,\xi _{\alpha },\eta ^{\alpha },g)$ be an trans-$S$%
-manifold and $\Omega $ \textit{the fundamental }$2$\textit{-form} of $%
M^{2n+s}$ defined by 
\begin{equation}
\Omega (X,Y)=g(X,fY),  \label{omega}
\end{equation}%
(see \cite{Nak-1966}). From Proposition 3.1 (i) in \cite{AFP-2016}, for a
trans-$S$-manifold, 
\begin{equation}
d\Omega =2\Omega \wedge \sum_{i=1}^{s}\beta _{i}\eta _{i}.  \label{domega}
\end{equation}%
If the fundamental $2$-form $\Omega $ on $M^{2n+s}$ is closed, then M
becomes a $K$-manifold. Moreover, $F=d\eta _{i}$, it becomes an $S-$%
manifold. If $d\eta _{i}=0$, it becomes a $C-$manifold. \ When $\Omega $ is
closed, the \textit{magnetic field} $F_{q}$ on $M^{2n+s}$ can be defined by 
\begin{equation*}
F_{q}(X,Y)=q\Omega (X,Y),
\end{equation*}%
where $X$ and $Y$ are vector fields on $M^{2n+s}$ and $q$ is a real
constant. $F_{q}$ is called the \textit{contact magnetic field with strength}
$q$ \cite{JMN-2015}. If $q=0$ then the magnetic curves are geodesics of $%
M^{2n+s}$. Because of this reason one can consider $q\neq 0$ (see \cite%
{CFG-2009} and \cite{DIMN-2015}).

From (\ref{F}) and (\ref{omega}), the Lorentz force $\Phi $ associated to
the contact magnetic field $F_{q}$ can be written as 
\begin{equation*}
\Phi _{q}=-qf.
\end{equation*}%
So the Lorentz equation (\ref{Lorentz eq}) can be written as%
\begin{equation}
\nabla _{T}T=-qfT,  \label{magneticcurve}
\end{equation}%
where $\gamma :I\subseteq R\rightarrow M^{2n+s}$ is a smooth unit-speed
curve and $T=\gamma ^{\prime }$ (see \cite{DIMN-2015} and \cite{JMN-2015}).

From \ref{domega}, for trans-$S$-manifolds, notice that $\Omega $ does not
need to be closed in general. But, we can still look for curves satisfying $%
\nabla _{T}T=-q\varphi T$ in a trans-$S$-manifold, calling them \textit{%
trajectories}. In this paper, for sake of computations, Legendre
trajectories will be considered. The general solution of the problem is in
progress.

For the last part of this study, it is necessary to define $C-$parallel $C-$%
proper curves as below:

We can generalize the definition from \cite{GO-2019} to trans-$S$-manifolds:

\begin{definition}
\cite{GO-2019}Let $\gamma :I\rightarrow (M^{2n+s},f,\xi _{i},\eta _{i},g)$
be a unit speed curve in an trans-$S$-manifold. Then $\gamma $ is called

i) $C$\textit{-parallel (in the tangent bundle) }if 
\begin{equation*}
\nabla _{T}H=\lambda \overset{s}{\underset{i=1}{\sum }}\xi _{i},
\end{equation*}

ii) $C$\textit{-parallel in the normal bundle }if 
\begin{equation*}
\nabla _{T}^{\perp }H=\lambda \overset{s}{\underset{i=1}{\sum }}\xi _{i},
\end{equation*}

iii) $C$\textit{-proper (in the tangent bundle) }if 
\begin{equation*}
\Delta H=\lambda \overset{s}{\underset{i=1}{\sum }}\xi _{i},
\end{equation*}

iv) $C$\textit{-proper in the normal bundle} if 
\begin{equation*}
\Delta ^{\perp }H=\lambda \overset{s}{\underset{i=1}{\sum }}\xi _{i},
\end{equation*}%
where $H$ is the mean curvature field of $\gamma $, $\lambda $ is a
real-valued non-zero differentiable function, $\nabla $ is the Levi-Civita
connection, $\nabla ^{\perp }$ is the Levi-Civita connection in the normal
bundle, $\Delta $ is the Laplacian and $\Delta ^{\perp }$ is the Laplacian
in the normal bundle.
\end{definition}

From the definition, same direct proposition as in \cite{GO-2019} is
obtained:

\begin{proposition}
\label{lemma1}\cite{GO-2019}Let $\gamma :I\rightarrow (M^{2n+s},f,\xi
_{i},\eta _{i},g)$ be a unit speed curve in an $S$-manifold. Then

i) $\gamma $ is $C$-parallel (in the tangent bundle) if and only if 
\begin{equation}
-\kappa _{1}^{2}T+\kappa _{1}^{\prime }E_{2}+\kappa _{1}\kappa
_{2}E_{3}=\lambda \overset{s}{\underset{i=1}{\sum }}\xi _{i},
\label{cparalleltangent}
\end{equation}

ii) $\gamma $ is $C$\textit{-parallel in the normal bundle }if and only if%
\begin{equation}
\kappa _{1}^{\prime }E_{2}+\kappa _{1}\kappa _{2}E_{3}=\lambda \overset{s}{%
\underset{i=1}{\sum }}\xi _{i},  \label{cparallelnormal}
\end{equation}

iii) $\gamma $ is $C$\textit{-proper (in the tangent bundle) }if and only if 
\begin{equation}
3\kappa _{1}\kappa _{1}^{\prime }T+\left( \kappa _{1}^{3}+\kappa _{1}\kappa
_{2}^{2}-\kappa _{1}^{\prime \prime }\right) E_{2}-(2\kappa _{1}^{\prime
}\kappa _{2}+\kappa _{1}\kappa _{2}^{\prime })E_{3}-\kappa _{1}\kappa
_{2}\kappa _{3}E_{4}=\lambda \overset{s}{\underset{i=1}{\sum }}\xi _{i},
\label{cpropertangent}
\end{equation}

iv) $\gamma $ is $C$\textit{-proper in the normal bundle} if and only if 
\begin{equation}
\left( \kappa _{1}\kappa _{2}^{2}-\kappa _{1}^{\prime \prime }\right)
E_{2}-\left( 2\kappa _{1}^{\prime }\kappa _{2}+\kappa _{1}\kappa
_{2}^{\prime }\right) E_{3}-\kappa _{1}\kappa _{2}\kappa _{3}E_{4}=\lambda 
\overset{s}{\underset{i=1}{\sum }}\xi _{i}.  \label{cpropernormal}
\end{equation}
\end{proposition}

\section{Main results on Legendre Trajectories}

Let $M=(M,f,\xi _{i},\eta _{i},g)$ be a trans-$S$-manifold and $\gamma
:I\rightarrow M$ a unit-speed Legendre curve with arc-length parameter $t$.
Assume that $\gamma $ satisfies $\nabla _{T}T=-qfT.$ Then, we have%
\begin{equation*}
\nabla _{T}T=-qfT=\kappa _{1}E_{2}
\end{equation*}%
and%
\begin{equation*}
g(fT,fT)=1.
\end{equation*}%
So, 
\begin{equation*}
fT\neq 0.
\end{equation*}%
Using the norm of both sides gives us%
\begin{equation}
\kappa _{1}=\left\vert q\right\vert .  \label{k1}
\end{equation}%
Thus%
\begin{equation*}
\left\vert q\right\vert E_{2}=-qfT
\end{equation*}%
and%
\begin{equation}
fT=\delta E_{2},  \label{fT}
\end{equation}%
where $\delta =sgn(-q)$. From (\ref{kappa1etae2}) and (\ref{fT}), we have%
\begin{equation*}
\beta _{i}\left\vert _{\gamma }\right. =0.
\end{equation*}%
(\ref{fT}) gives us%
\begin{equation}
fE_{2}=-\delta T.  \label{fE2}
\end{equation}%
From (\ref{nablafT}) and (\ref{fE2}), we can write \ 
\begin{eqnarray*}
\nabla _{T}fT &=&\delta \nabla _{T}E_{2}=\delta (-\kappa _{1}T+\kappa
_{2}E_{3}) \\
&=&\sum_{i=1}^{s}\alpha _{i}\xi _{i}-\kappa _{1}\delta T.
\end{eqnarray*}%
As a result, we find%
\begin{equation}
\kappa _{2}E_{3}=\delta \sum_{i=1}^{s}\alpha _{i}\xi _{i},  \label{kappa2e3}
\end{equation}%
which gives us%
\begin{equation}
\kappa _{2}=\sqrt{\sum_{i=1}^{s}\alpha _{i}^{2}}.  \label{kappa2}
\end{equation}%
Then%
\begin{equation*}
\kappa _{2}=0\Leftrightarrow \alpha _{i}\left\vert _{\gamma }\right. =0.
\end{equation*}%
Let $\kappa _{2}\neq 0$. Notice that $sgn\left( g\left(
E_{3},\sum_{i=1}^{s}\alpha _{i}\xi _{i}\right) \right) =\delta $. Using (\ref%
{kappa2e3}) and (\ref{kappa2}), we find%
\begin{equation}
E_{3}=\frac{\delta }{\sqrt{\sum_{i=1}^{s}\alpha _{i}^{2}}}%
\sum_{i=1}^{s}\alpha _{i}\xi _{i}.  \label{E3}
\end{equation}%
If we differentiate $E_{3}$, we obtain%
\begin{equation}
\kappa _{3}E_{4}=\delta \sum_{i=1}^{s}(\frac{\alpha _{i}}{\sqrt{%
\sum_{i=1}^{s}\alpha _{i}^{2}}})^{\prime }\xi _{i}.  \label{E4}
\end{equation}%
\begin{equation}
\kappa _{3}=\sqrt{\sum_{i=1}^{s}[(\frac{\alpha _{i}}{\sqrt{%
\sum_{i=1}^{s}\alpha _{i}^{2}}})^{\prime }]^{2}}.  \label{K3}
\end{equation}%
Moreover, if $\kappa _{3}=0$, then%
\begin{equation*}
\frac{\alpha _{i}}{\sqrt{\sum_{i=1}^{s}\alpha _{i}^{2}}}=c_{i}=constant\text{%
, }\forall i.
\end{equation*}%
Hence%
\begin{equation*}
\sum_{i=1}^{s}\alpha _{i}^{2}\left( \sum_{i=1}^{s}c_{i}^{2}-1\right) =0.
\end{equation*}%
So, 
\begin{equation*}
\sum_{i=1}^{s}\alpha _{i}^{2}=0\Leftrightarrow \kappa _{2}=0,
\end{equation*}%
or%
\begin{equation*}
\sum_{i=1}^{s}c_{i}^{2}=1.
\end{equation*}%
To sum up, if $\kappa _{3}=0$ and $\kappa _{2}\neq 0$, we have%
\begin{equation*}
E_{2}=\delta fT,
\end{equation*}%
\begin{equation*}
E_{3}=\delta \sum_{i=1}^{s}c_{i}\xi _{i},
\end{equation*}%
where 
\begin{equation}
\alpha _{i}=c_{i}\sum_{i=1}^{s}\alpha _{i}^{2},\forall i,  \label{equ1}
\end{equation}%
\begin{equation}
c_{i}=\text{constant such that }\sum_{i=1}^{s}c_{i}^{2}=1.  \label{equ2}
\end{equation}%
Now we can state the following theorem:

\begin{theorem}
Let $\gamma :I\rightarrow M$ be a Legendre trajectory. Then $\gamma $ is one
of the following:

$1)$ a Legendre circle with $\kappa _{1}=\left\vert q\right\vert $ and the
Frenet frame field $\left\{ T,\delta fT\right\} ,$ where $\delta =sgn(-q)$.
In this case, $\alpha _{i}=0,$ $\beta _{i}=0,\forall i.$

$2)$ a Legendre curve of osculating order $r\geq 3$ with 
\begin{equation*}
\kappa _{1}=\left\vert q\right\vert ,\kappa _{2}=\sqrt{\sum_{i=1}^{s}\alpha
_{i}^{2}},
\end{equation*}
$\kappa _{3}$ given in $(\ref{K3})$ and the Frenet frame field 
\begin{equation*}
\left\{ T,\delta fT,E_{3},E_{4},...,E_{r}\right\} ,
\end{equation*}
where $\delta =sgn(-q);$ $E_{3}$, $E_{4}$ are given in (\ref{E3}) and (\ref%
{E4}), respectively. In this case, $\alpha _{i}\neq 0,$ $\exists i,$ $\beta
_{i}=0,\forall i.$ Moreover, if $r=3$, equations $(\ref{equ1})$ and $(\ref%
{equ2})$ are also satisfied and its Frenet frame field is 
\begin{equation*}
\left\{ T,\delta fT,E_{3}\right\} .
\end{equation*}
\end{theorem}

\section{Main results of $C$-parallel and $C$-proper Legendre Curves}

Let $M^{2n+s}$ be $a$ trans-S-manifold and $\gamma :I\rightarrow M$ a
Legendre curve in $M$.

\textbf{i) }$C$-parallel in the tangent bundle:%
\begin{equation*}
-\kappa _{1}^{2}T+\kappa _{1}^{\prime }E_{2}+\kappa _{1}\kappa
_{2}E_{3}=\lambda \overset{s}{\underset{i=1}{\sum }}\xi _{i}.
\end{equation*}%
If we apply $T$ to both sides, we have the following result:

\begin{theorem}
There does not exist a $C$-parallel Legendre curve (in the tangent bundle)
in a trans-$S$-manifold.
\end{theorem}

\textbf{ii)} $C$-parallel in the normal bundle:%
\begin{equation*}
\kappa _{1}^{\prime }E_{2}+\kappa _{1}\kappa _{2}E_{3}=\lambda \overset{s}{%
\underset{i=1}{\sum }}\xi _{i}.
\end{equation*}

\textbf{a)} $r=2$.

\begin{equation*}
\kappa _{1}^{\prime }E_{2}=\lambda \overset{s}{\underset{i=1}{\sum }}\xi
_{i}.
\end{equation*}

\begin{theorem}
Let $r=2$. Then $\gamma $ is $C$-parallel in the normal bundle iff%
\begin{equation*}
\kappa _{1}=\mp \sqrt{s}\beta ,
\end{equation*}%
\begin{equation*}
\lambda =-\beta ^{\prime },
\end{equation*}%
\begin{equation*}
\overset{s}{\underset{i=1}{\sum }}\xi _{i}=\pm \sqrt{s}E_{2}.
\end{equation*}%
In this case, $\beta _{1}=\beta _{2}=...=\beta _{s}=\beta .$
\end{theorem}

\textbf{b)} $r\geq 3$.

In this case, for a smooth function $w=w(t)$, we have%
\begin{equation}
\overset{s}{\underset{i=1}{\sum }}\xi _{i}=\sqrt{s}\left( \cos wE_{2}+\sin
wE_{3}\right) .  \label{equa1}
\end{equation}%
If we differentiate the above equation, we have%
\begin{equation}
\overset{s}{\underset{i=1}{\sum }}\beta _{i}=-\sqrt{s}\kappa _{1}\cos w
\label{equa2}
\end{equation}%
and%
\begin{equation}
\kappa _{2}=\pm \frac{1}{\sqrt{s}}\overset{s}{\underset{i=1}{\sum }}\alpha
_{i}-w^{\prime }.  \label{equa3}
\end{equation}%
We also have%
\begin{equation}
\lambda =\frac{-\kappa _{1}\kappa _{1}^{\prime }}{\overset{s}{\underset{i=1}{%
\sum }}\beta _{i}}.  \label{equa4}
\end{equation}%
Since $fT\perp E_{2}$, we can write%
\begin{equation}
fT=\pm \left( \sin wE_{2}-\cos wE_{3}\right) .  \label{equa5}
\end{equation}

\begin{theorem}
Let $r\geq 3$. Then $\gamma $ is $C$-parallel in the normal bundle iff
equations $(\ref{equa1})$,$(\ref{equa2})$,$(\ref{equa3})$,$(\ref{equa4})$
and $(\ref{equa5})$ are satisfied.
\end{theorem}

\textbf{iii)} $C$-proper in the tangent bundle:%
\begin{equation*}
3\kappa _{1}\kappa _{1}^{\prime }T+\left( \kappa _{1}^{3}+\kappa _{1}\kappa
_{2}^{2}-\kappa _{1}^{\prime \prime }\right) E_{2}-(2\kappa _{1}^{\prime
}\kappa _{2}+\kappa _{1}\kappa _{2}^{\prime })E_{3}-\kappa _{1}\kappa
_{2}\kappa _{3}E_{4}=\lambda \overset{s}{\underset{i=1}{\sum }}\xi _{i}.
\end{equation*}%
If we apply $T$, we directly have $\kappa _{1}=$constant. Then the equation
reduces to%
\begin{equation*}
\kappa _{1}\left( \kappa _{1}^{2}+\kappa _{2}^{2}\right) E_{2}-\kappa
_{1}\kappa _{2}^{\prime }E_{3}-\kappa _{1}\kappa _{2}\kappa
_{3}E_{4}=\lambda \overset{s}{\underset{i=1}{\sum }}\xi _{i}.
\end{equation*}%
Applying $E_{2},$ we get%
\begin{equation*}
\kappa _{1}^{2}\left( \kappa _{1}^{2}+\kappa _{2}^{2}\right) =-\lambda 
\overset{s}{\underset{i=1}{\sum }}\beta _{i}.
\end{equation*}

\textbf{a)} $r=2$.

In this case, we have%
\begin{equation}
\kappa _{1}^{3}E_{2}=\lambda \overset{s}{\underset{i=1}{\sum }}\xi _{i}.
\label{eq1}
\end{equation}%
If we apply $\xi _{j}$, we find%
\begin{equation*}
\kappa _{1}^{3}\eta _{j}\left( E_{2}\right) =\lambda ,\forall j.
\end{equation*}%
If we denote%
\begin{equation*}
\beta _{1}=\beta _{2}=...=\beta _{s}=\beta ,
\end{equation*}%
we get%
\begin{equation*}
\lambda =-s\beta ^{3}=constant.
\end{equation*}%
If we differentiate (\ref{eq1}), it is easy to see that%
\begin{equation*}
\overset{s}{\underset{i=1}{\sum }}\alpha _{i}=0.
\end{equation*}%
As a result, we have%
\begin{equation*}
\kappa _{1}=\mp \sqrt{s}\beta =\text{constant,}
\end{equation*}%
i.e., $\gamma $ is a circle.

\begin{theorem}
Let $r=2$. Then $\gamma $ is $C$-proper in the tangent bundle iff it is a
circle with 
\begin{equation*}
\kappa _{1}=\mp \sqrt{s}\beta =constant
\end{equation*}
and the Frenet frame field 
\begin{equation*}
\left\{ T,\frac{\pm 1}{\sqrt{s}}\overset{s}{\underset{i=1}{\sum }}\xi
_{i}\right\} .
\end{equation*}
In this case, $\beta _{1}=\beta _{2}=...=\beta _{s}=\beta ,$ $\lambda
=-s\beta ^{3}=constant$ and $\overset{s}{\underset{i=1}{\sum }}\alpha
_{i}=0. $
\end{theorem}

\bigskip

\textbf{b)} $r=3$.%
\begin{equation*}
\kappa _{1}\left( \kappa _{1}^{2}+\kappa _{2}^{2}\right) E_{2}-\kappa
_{1}\kappa _{2}^{\prime }E_{3}=\lambda \overset{s}{\underset{i=1}{\sum }}\xi
_{i}.
\end{equation*}%
So, $\overset{s}{\underset{i=1}{\sum }}\xi _{i}\in sp\left\{
E_{2},E_{3}\right\} $. It can be written as%
\begin{equation}
\overset{s}{\underset{i=1}{\sum }}\xi _{i}=\sqrt{s}\left( \cos wE_{2}+\sin
wE_{3}\right) ,  \label{equat1}
\end{equation}%
for a smooth function $w=w(t)$. If we differentiate this equation and apply $%
T$, we find 
\begin{equation}
\overset{s}{\underset{i=1}{\sum }}\beta _{i}=-\sqrt{s}\kappa _{1}\cos w,
\label{equat2}
\end{equation}%
and%
\begin{equation}
\kappa _{2}=\pm \frac{1}{\sqrt{s}}\overset{s}{\underset{i=1}{\sum }}\alpha
_{i}-w^{\prime }.  \label{equat3}
\end{equation}%
We also have%
\begin{equation}
fT=\pm \left( \sin wE_{2}-\cos wE_{3}\right)  \label{equat4}
\end{equation}%
and%
\begin{equation}
\lambda =\frac{-\kappa _{1}^{2}\left( \kappa _{1}^{2}+\kappa _{2}^{2}\right) 
}{\overset{s}{\underset{i=1}{\sum }}\beta _{i}}.  \label{equat5}
\end{equation}

\begin{theorem}
Let $r=3$. Then $\gamma $ is $C$-proper in the tangent bundle iff equations $%
\left( \ref{equat1}\right) $, $\left( \ref{equat2}\right) $, $\left( \ref%
{equat3}\right) $, $\left( \ref{equat4}\right) $ and $\left( \ref{equat5}%
\right) $ are satisfied.
\end{theorem}

\textbf{c) }$\mathbf{r}\geq 4$.

In this case, $\overset{s}{\underset{i=1}{\sum }}\xi _{i}\in sp\left\{
E_{2},E_{3},E_{4}\right\} $, consequently, $fT\in sp\left\{
E_{2},E_{3},E_{4},E_{5}\right\} .$ Let us write%
\begin{equation}
\overset{s}{\underset{i=1}{\sum }}\xi _{i}=\sqrt{s}\left( \cos wE_{2}+\sin
w\cos \varphi E_{3}+\sin w\sin \varphi E_{4}\right)  \label{eq2}
\end{equation}%
for some smooth functions $w=w(t)$, $\varphi =\varphi (t).$ As a result, the
curve must satisfy%
\begin{equation*}
\kappa _{1}=constant,
\end{equation*}%
\begin{equation*}
\overset{s}{\underset{i=1}{\sum }}\beta _{i}=-\sqrt{s}\kappa _{1}\cos w,
\end{equation*}%
\begin{equation*}
\lambda \sqrt{s}\cos w=\kappa _{1}^{2}\left( \kappa _{1}^{2}+\kappa
_{2}^{2}\right) ,
\end{equation*}%
\begin{equation*}
\lambda \sqrt{s}\sin w\cos \varphi =-\kappa _{1}\kappa _{2}^{\prime },
\end{equation*}%
\begin{equation*}
\lambda \sqrt{s}\sin w\sin \varphi =-\kappa _{1}\kappa _{2}\kappa _{3}.
\end{equation*}%
Differentiating (\ref{eq2}), we also have%
\begin{equation*}
\kappa _{4}=\frac{-\left( \overset{s}{\underset{i=1}{\sum }}\alpha
_{i}\right) .g\left( fT,E_{5}\right) }{\sin w\sin \varphi }.
\end{equation*}

\begin{theorem}
Let $r\geq 4$. Then $\gamma $ is $C$-proper in the tangent bundle iff it
satisfies the last six equations.
\end{theorem}

\textbf{iv)} $C$-proper in the normal bundle:

\begin{equation*}
\left( \kappa _{1}\kappa _{2}^{2}-\kappa _{1}^{\prime \prime }\right)
E_{2}-\left( 2\kappa _{1}^{\prime }\kappa _{2}+\kappa _{1}\kappa
_{2}^{\prime }\right) E_{3}-\kappa _{1}\kappa _{2}\kappa _{3}E_{4}=\lambda 
\overset{s}{\underset{i=1}{\sum }}\xi _{i}.
\end{equation*}%
In this case, again $\overset{s}{\underset{i=1}{\sum }}\xi _{i}\in sp\left\{
E_{2}\right\} ,$ $\overset{s}{\underset{i=1}{\sum }}\xi _{i}\in sp\left\{
E_{2},E_{3}\right\} $ or $\overset{s}{\underset{i=1}{\sum }}\xi _{i}\in
sp\left\{ E_{2},E_{3},E_{4}\right\} $ depending on the osculating order $r$.
We can follow the above procedure to get results for $r=2\ $and $r=3$. The
case $r\geq 4$ is similar to case \textbf{iii) c)} with minor changes in
equations.

\bigskip

\textbf{Remark. }For sake of shortness\textbf{, }$\alpha _{i}\left\vert
_{\gamma }\right. $ and $\beta _{i}\left\vert _{\gamma }\right. $ are
written as $\alpha _{i}$ and $\beta _{i}$ where possible. This means the
equations are not necessarily satisfied globally. But instead, they are
satisfied along the curve $\gamma .$

\end{document}